\documentclass[12pt,a4paper]{article}

\usepackage[latin1]{inputenc}
\usepackage{color}

\usepackage{bussproofs}
\usepackage{graphicx}
\usepackage{graphics}
\usepackage{latexsym}
\usepackage{amssymb,amsmath}
\usepackage{euscript}
\usepackage{enumerate}
\usepackage{amsmath}
\usepackage{amsfonts}
\usepackage{multicol}
\usepackage{etex}
\usepackage[all]{xy}
\usepackage{color}

\usepackage[colorlinks=true, pdfstartview=FitV, linkcolor=blue, citecolor=blue, urlcolor=blue, 
]
{hyperref}
\newtheorem{teo}{Theorem}[section]
\newtheorem{coro}[teo]{Corollary}
\newtheorem{lema}[teo]{Lemma}

\newtheorem{obs}[teo]{Observation}

\newtheorem{conj}[teo]{Conjecture}
\newtheorem{prob}[teo]{Problem}
\newenvironment{dem}{\noindent {\em Proof.} }{\hfill $\square$ \bigskip}
\newenvironment{prueba}[1]{\noindent\bf Proof of #1. \rm}{$\quad \hfill \square$ \bigskip}

\newcommand{\R}{\mathbb{R}}

\newcommand{\Z}{\mathbb{Z}}

\newcommand{\essinf}{\mathop{\text{ess inf}}}

\def\XXint#1#2#3{{\setbox0=\hbox{$#1{#2#3}{\int}$}
    \vcenter{\hbox{$#2#3$}}\kern-.5\wd0}}

\usepackage{fancyhdr}

\fancypagestyle{encabezadotesis}{
\fancyhead{}
\fancyhead[RO]{\em \nouppercase{\rightmark}}
\fancyhead[LE]{\bf \nouppercase{\leftmark}}
\cfoot{\thepage}
}

\begin{document}
\title{Quantitative weighted mixed weak-type inequalities for classical operators}
\author{S. Ombrosi, C. P\'erez and J. Recchi}
\date{}
\maketitle

\renewcommand{\thefootnote}{\fnsymbol{footnote}}
 \footnotetext{
The first and  third authors are  supported by Universidad Nacional del Sur and CONICET.}
\footnotetext{The second author was
supported by grant MTM2014-53850-P, Spanish Government}

\begin{abstract}

 We improve on several mixed weak type inequalities
 both for the Hardy-Littlewood maximal function and for Calder\'on-Zygmund operators.
  These type of inequalities were considered by Muckenhoupt and Wheeden and later on
 by Sawyer estimating the $L^{1, \infty}(uv)$ norm of $v^{-1}T(fv)$ for special cases. The emphasis is made in proving new and more precise quantitative estimates involving the $A_p$ or $A_\infty$ constants of the weights involved.

\end{abstract}

\section{Introduction and statements of the main results}

Let $M$ denote the usual Hardy-Littlewood maximal function, then  according to a fundamental result of B. Muckenhoupt  \cite{M 72}, $M$ is a bounded operator on the Lebesgue space $L^p(d\mu)$, $1 <p < +\infty$, if and only if $d\mu = w(x)dx$ and the weight $w$ satisfies the simple geometric condition
$$
[w]_{A_p} :=  \sup_Q \left(\frac{1}{|Q|}\int_Q w\right)\left(\frac{1}{|Q|}\int_Q w^{1-p'}\right)^{p-1} <\infty,
$$
where the supremum is taken over all cubes $Q$ in $\R^n$. This is the celebrated Muckenhoupt $A_p$ condition.  A similar result holds in the case $p = 1$, namely $M$ is of weak type
(1,1) with respect to $\mu$, i.e. $
M:L^1(\mu)\to L^{1,\infty}(\mu)$, if and only if $d\mu = w(x)dx$ and the weight $w$ satisfies the
$A_1$ condition,
$$
[w]_{A_1}:=  \sup_Q \left(\frac{1}{|Q|}\int_Q w\right)(\essinf_Q w)^{-1} <\infty
% \esssup_{x \in \R^n} \frac{Mw(x)}{w(x)}
$$
where, again,  the supremum is taken over all cubes $Q$ in $\R^n$.

Since the $A_p$ theorem of Muckenhoupt plays a central role in modern Harmonic Analysis,  different proofs from the original one in \cite{M 72} have been considered in the literature.
In particular, E. Sawyer tried in \cite{Sawyer85} the following approach based on the factorization theorem for $A_p$ weights of P. Jones (see \cite{garcia curva-rubio de francia}).  Recall that a weight $w$ satisfies the $A_p$ condition if and only if there are two $A_1 $ weights $u$ and $v$ such that
\begin{equation}\label{desigualdad 2 de sawyer}
w=uv^{1-p}.
\end{equation}
Then, if the following operator is defined
$$Sf=\frac{M(vf)}{v}$$
the boundedness of $M$ on $L^p(w)$ may be rewritten as
\begin{equation}\label{desigualdad 3 de sawyer}
\int_{\R^n} |Sf|^p \,uvdx\leq c \int_{\R^n} |f|^p \,uvdx.
\end{equation}
Observe now that since $v\in A_1$, $Mv \leq [v]_{A_1} v$ and hence $S$ is bounded in $L^{\infty}(uv)$. Therefore, if we show that $S$ is of weak type $(1,1)$ with respect to the measure $uvdx$ we can apply the Marcinkiewicz interpolation theorem to derive  \eqref{desigualdad 3 de sawyer}.
This is precisely the statement of the following theorem from \cite{Sawyer85}.
%This approach was considered by Sawyer  in \cite{Sawyer85} where it is proved the following result.

\begin{teo}\label{teorema de sawyer}
If  $u$, $v \in A_1(\R)$, then
$$
\left\|\frac{M(g)}{v}\right\|_{L^{1, \infty}(uv)} \leq c\;\|g\|_{L^1(u)},
$$
where $c$ depends only on the $A_1$ constant of $u$ and the $A_1$ constant of $v$.   This shows that the operator  $Sf=v^{-1}M(vf)$ is of weak type  $(1,1)$ with respect to the measure  $vudx$.
\end{teo}

In the same article, Sawyer conjectured that this theorem should also hold for the maximal function in $\R^n$  and for the the Hilbert transform $H$ instead of $M$.

The article of Sawyer was also very much motivated by a previous work of B. Muckenhoupt and R. Wheeden  \cite{M-W 77}. The main result of this paper holds this time for both the one dimensional Hardy-Littlewood maximal function and the Hilbert transform. To be more precise, the  main result proved in \cite{M-W 77} is the following.

\color{black}

\begin{teo}\label{teorema 1.1 de cump}
Let $w\in A_1(\R)$, there exists a constant  $c$ such that,
\begin{equation}
\left\|  M(fw^{-1})w \right\|_{L^{1, \infty}(\R)}
\leq c\, \left\|  f\right\|_{L^{1}(\R)}\,
\end{equation}
and
\begin{equation}
\left\|  H(fw^{-1})w \right\|_{L^{1, \infty}(\R)}  \leq c\, \left\|  f\right\|_{L^{1}(\R)}.
\end{equation}
\end{teo}
%

%\comment{CARLOS  Acá agrego la definición de $A_\infty(u)$ }

In \cite{cruz-uribe martell perez} the authors extended both Theorems \ref{teorema de sawyer}  and \ref{teorema 1.1 de cump} 
to $\R^n$ containing in particular the conjectures formulated by Sawyer mentioned above. The precise result is the following.
%More precisely, they proved the following theorem.

%
\begin{teo}\label{teorema 1.3 de cump}
Suppose that $u \in A_1$ and that either  $v \in A_1$ or $v \in A_{\infty}(u)$, then there exists a constant $c$ such that,
\begin{equation}\label{sawyer en mas dimensiones}
\left\|\frac{M(fv)}{v}\right\|_{L^{1, \infty}(uv)} \leq c\;\|f\|_{L^1(uv)}
\end{equation}
and
\begin{equation} \label{sawyer para T}
\left\|\frac{T(fv)}{v}\right\|_{L^{1, \infty}(uv)}\leq c\;\|f\|_{L^1(uv)},
\end{equation}
where $M$ is the Hardy-Littlewood maximal operator and $T$ is a Calderón-Zygmund operator.
\end{teo}

We remark that this result holds for $T^*$, the maximal singular integral operator, instead of $T$. Given weights $u$ and $v$, by $v \in A_\infty (u)$ we mean that $v$ satisfies the $A_\infty$ condition defined with respect to the measure $u dx$ (as opposed to Lebesgue measure). A more precise definition is given in Section \ref{preliminares} below.

We emphasize that this theorem contains both Theorems \ref{teorema de sawyer} and  \ref{teorema 1.1 de cump} as particular cases. 
 Indeed,  the case of the first theorem is clear. For the second, if $w \in A_1$, we let $u=w$ and $v=w^{-1}$. Then,  $uv = 1 \in A_\infty$ and thus $v \in A_\infty(u)$ by Lemma \ref{lema 2.1 de CUMP} and Observation \ref{observacion 2.2 de CUMP}.

To prove Theorem \ref{teorema 1.3 de cump}, the authors show that it suffices to prove the result for the dyadic maximal function $M_d$ by proving an extrapolation type theorem, Theorem \ref{teorema 1.7 de cump} below, that allows to replace $T$ or $M$ by $M_d$. To be more precise, the combination of the following  two theorems from \cite{cruz-uribe martell perez} proves Theorem \ref{teorema 1.3 de cump}.

\begin{teo}\label{teorema 1.4 de cump}
Suppose that $u \in A_1$ and that either  $v \in A_1$ or $v \in A_{\infty}(u)$, then there exists a constant $c$ such that,
\begin{equation}\label{sawyer en mas dimensiones para diadica}
\left\|\frac{M_d(fv)}{v}\right\|_{L^{1, \infty}(uv)} \leq c\;\|f\|_{L^1(uv)}.
\end{equation}
\end{teo}
%Para probar la desigualdad (\ref{sawyer en mas dimensiones para diadica}) cuando $u \in A_1$ y $v\in A_\infty(u)$, ellos usan una simple descomposición de Calderón-Zygmund inspirada en la prueba que dan Muckenhoupt y Wheeden del Teorema \ref{teorema 1.1 de cump}.
%Para el caso $u\in A_1$ y $v\in A_1$, la prueba es una adaptación de la dada por Sawyer en la recta, la cual se realiza usando un muy delicado argumento de descomposición.

%El segundo paso para probar el Teorema \ref{teorema 1.3 de cump} fue obtener un teorema de tipo extrapolación usando técnicas  de \cite{cruz-uribe martell perez 04}. El teorema es el siguiente,

\begin{teo}\label{teorema 1.7 de cump}
Given a family $\mathcal{F}$ of pair of functions, suppose that for some  $p\in (0,\infty)$  and for every $w\in A_\infty$,
$$
\|f\|_{L^p(w)}\leq C \|g\|_{L^p(w)},
$$
for all $(f,g) \in \mathcal{F} $ such that the left-hand side is finite, and where $C$ depends only on the $A_\infty$ constant of $w$. Then for all weights   $u \in A_1$ and $v \in A_\infty$,
$$
\|fv^{-1}\|_{L^{1,\infty}(uv)}\leq C \|gv^{-1}\|_{L^{1,\infty}(uv)}\;\;\;\; (f,g) \in \mathcal{F}.
$$
Here $\mathcal{F}$ denotes a family of ordered pairs of non-negative, measurable functions $(f,g)$.
\end{teo}

%The proof of Theorem \ref{teorema 1.3 de cump} is now immediate. 

Theorem \ref{teorema 1.7 de cump} from \cite{cruz-uribe martell perez} is used  to pass from $M$ to $M_d$ since by standard methods,  for every  $p\in (0,\infty)$ and every   $w \in A_\infty$
$$
||M(fv)||_{L^p(w)}\leq c\;||M_d(fv)||_{L^p(w)},
$$
where the constant $c$ involves the $A_\infty$ constant of  $w$. However, there are recent results showing that Theorem \ref{teorema 1.7 de cump} can be avoided in the transition from $M$ to $M_d$. 
%being the constant $C$ dimensional constant. 
Indeed, using for instance \cite{hytonen-perez} p. 792,  we have that 
\begin{equation*}
  Mf\leq c_n\sum_{\alpha\in\{0,\frac13\}^n}M_d^{\alpha}f.
\end{equation*}
where $M_d^{\alpha}$ is  an appropriate shifted dyadic maximal function with similar properties as $M_d$. Thus, the expression on the left in \eqref{sawyer en mas dimensiones} is bounded by a dimensional constant multiple of the corresponding expression for $M_d^{\alpha}$. Since each of these $M_d^{\alpha}$ has similar properties as $M_d$  the corresponding proof of \eqref{sawyer en mas dimensiones para diadica} is exactly the same.

In \cite{cruz-uribe martell perez} the authors conjectured that Theorem \ref{teorema 1.4 de cump} still holds under milder hypotheses on the weight $v$. To be more precise, the authors state what is now known as \lq\lq Sawyer's Conjecture\rq\rq, although E. Sawyer never asserted it. The conjecture is the following.

\begin{conj}\label{conjetura de sawyer}
Suppose that $u \in A_1$ and $v \in A_\infty$. Then there exists a constant $c$ such that
\begin{equation}
\left\| \frac{M_d(fv)}{v}\right\|_{L^{1,\infty}(uv)}\leq c\,\|f\|_{L^1(uv)}.
\end{equation}
\end{conj}

Note that if $v \in A_\infty(u)$ (always assuming $u\in A_1$), then $v \in A_\infty$ (see Lemmas \ref{lema 2.1 de CUMP} and \ref{lema 2.5 de cump}). This conjecture has been open for several years and has been studied by different authors.

In this paper we try to understand the difficulties of this conjecture and propose alternative ways to prove it. We will also study how the constants of the weights $u$ and $v$ is reflected in these inequalities, that is, we look for quantitative versions of this type of inequalities.

The first question that we pose concerning Sawyer's Theorem is the following:
\begin{quote}
What is the sharp dependence on the constants of the weights $u$ and $v$ when both are in $A_1$?
\end{quote}

Following the proof given in \cite{cruz-uribe martell perez},  which is an adaptation of the original proof given by Sawyer in \cite{Sawyer85} for the real line, we show  the dependence on the weight constants. More specifically, we prove the following result.

\begin{teo} \label{terorema de costantes cuadrado y a la cuarta}
If $u \in A_1$ and $v \in A_1$, there exists a dimensional constant $c$ such that
\begin{equation}\label{2-4estimate}
\left\| \frac{M_d (fv)}{v} \right\|_{L^{1, \infty}(uv)} \; \leq \; c\; [u]_{A_1}^2[v]_{A_1}^4\|f\|_{L^1(uv)}.
\end{equation}
\end{teo}
The proof may be found in Section \ref{dependencia de las constantes u y v}.

 We believe that the dependency on the constants in inequality \eqref{2-4estimate} is not sharp since the method does not seem to be adequate.  Trying to understand this issue we will focus on the special case $u = 1$ which is interesting in its own.  The finiteness of the estimate in this special case is assured by Theorem \ref{teorema 1.3 de cump} assuming even a weaker condition on $v$ than $A_1$, namely $v\in A_{\infty}(u)=A_{\infty}$. The method that we use is different from the one considered in the proof of Theorem \ref{terorema de costantes cuadrado y a la cuarta} allowing us to obtain more precise estimates. In particular we will prove the linearity of the constant bound of the weight $v$ if we assume the stronger condition $v\in A_1$ and the result is sharp. Our theorem is the following.

\begin{teo} \label{teo u y v en A1 para la recta}
Let $v \in A_1$. There exists a dimensional constant $c$, independent from $[v]_{A_1}$, such that

\begin{equation} \label{casevA1}
\left\|\frac{M(f)}{v}\right\|_{L^{1,\infty}(v)}\leq c \;[v]_{A_1} \|f\|_{L^1(\R^n)}.
\end{equation}
Furthermore, the linear dependence on \,$[v]_{A_1}$\,  is sharp.
\end{teo}

However, we want to understand the more general case.
\begin{prob} \label{teorema v en Ainfty}
Find an increasing function $\varphi: [1,\infty] \to [1,\infty]$ for which the following inequality holds whenever  $v \in A_\infty$
\begin{equation} \label{casevAinfty}
\left\| \frac{M (f)}{v} \right\|_{L^{1, \infty}(v)} \; \leq \; c\,\varphi([v]_{A_\infty}) \|f\|_{L^1(\R^n)},
\end{equation}
where $c$ is a constant that depends on the dimension.
\end{prob}

 This problem is a special case of Conjecture \ref{conjetura de sawyer} with $u=1$
and it will be studied In Section \ref{Caso especial de la conjetura de Sawyer para $M_d$}. 
The best constant in \eqref{casevAinfty}, \,$\varphi([v]_{A_\infty})$,\, is finite by Theorem \ref{teorema 1.3 de cump}. 
Our goal is to determine the best dependence on the constant of the weight $v$ or, in other words, to find the smallest function $\varphi$. Recall that \, $A_\infty = \cup_{p\geq1}A_p$ and that, if $w\in A_\infty$ we use the weight constant
\begin{equation} \label{constantAinfty}
[w]_{A_\infty}:= \sup_Q \frac{1}{w(Q)}\int_Q M(\chi_Q w)\,dx,
\end{equation}
called the Fujii-Wilson constant in some recent papers. We could use instead the constant defined by Hrushev in \cite{Hruscev}
 which is more natural, however it was shown in \cite{hytonen-perez} that it is much larger than the one given  by the functional \eqref{constantAinfty}.

We remark here that a condition on the weight $v$ in \eqref{casevA1} or \eqref{casevAinfty} must be taken into account. Indeed, there are estimates like
\begin{equation} \label{casev=Mw}
\left\|\frac{M(f)}{Mw}\right\|_{L^{1,\infty}(Mw)}\leq c\, \|f\|_{L^1(\R^n)},
\end{equation}
namely with $v=Mw$, which are false for a general function $w$ or measure.  This will be shown in Section
\ref{contraejemplo} where, furthermore,   an interesting relationship with the two weight problem for singular integrals is implicit in the argument.  In general, {weights of the form $Mw$} are not $A_{\infty}$ weights but small perturbations, namely when $v=(Mw)^{\delta}$, $\delta \in(0,1)$,  makes the inequality to be true since in this case $v\in A_1$ and Theorem \ref{casevA1} applies. It is interesting that
in special situations and for large perturbations of the weight the result is still true.  Indeed, if $v(x)=|x|^{-nr}  \approx (M\delta)^{r}$ with $r>1$, then there is a finite constant $c$ such that
\begin{equation} \label{max1}
\left\|\frac{M(f)}{v}\right\|_{L^{1,\infty}(v)}\leq c\, \|f\|_{L^1(\R^n)},
\end{equation}
being the result false in the case $r=1$.   This was proved in dimension one by Andersen and Muckenhoupt in \cite{AM} and by Mart{\'\i}n-Reyes, Ortega Salvador and
Sarri{\'o}n Gavi{\'a}n \cite{MOS} in higher dimensions. We remark that these weights $v(x)=|x|^{-nr}$ are not $A_{\infty}$ weights.

In view of Theorem  \ref{teo u y v en A1 para la recta} and the case $v=1$ we state the following conjecture for the general case.
%return now to Sawyer's Theorem for $u \in A_1$ y $v \in A_1$ and look for the sharp dependence on the weight constants. Our conjecture is the following.

\begin{conj}
Let $u \in A_1$ and $v\in A_1$, then there exists a dimensional constant $c$ such that
$$
\left\| \frac{M_d (fv)}{v} \right\|_{L^{1, \infty}(uv)}  \; \leq \; c\; [u]_{A_1}[v]_{A_1}\|f\|_{L^1(uv)}.
$$
\end{conj}

To see that the dependency cannot be better than $[u]_{A_1}[v]_{A_1}$ we prove the following result which strengthens our conjecture.

\begin{teo}\label{teo que las constantes son como mimimo el producto}
Let $u \in A_1, v \in A_1$. If
$$
\left\| \frac{M_d (fv)}{v} \right\|_{L^{1, \infty}(uv)}\leq c \,\varphi([u]_{A_1},[v]_{A_1})\|f\|_{L^1(uv)},
$$
then, there is a constant $c$ independent of the weights such that
$$
\varphi([u]_{A_1},[v]_{A_1})\geq c\,[u]_{A_1}[v]_{A_1}.
$$
\end{teo}

Another related problem, partly  intermediate  between the previous two problems,  would be to determine how the dependence on the constant $[v]_{A_p}$ is if we assume that $v \in A_p$ for some $p \geq 1$. We should also take into account that Theorem \ref{teo u y v en A1 para la recta} gives the sharp dependence on the real line when assuming the stronger assumption  $v \in A_1$. Based on this we state  the following conjecture.
\begin{conj} Let $v \in A_p$, $p \geq 1$, then there exists a dimensional constant $c$ such that
$$
\left\| \frac{M_d (fv)}{v} \right\|_{L^{1, \infty}(v)}  \; \leq \; c \,[v]_{A_p} \|f\|_{L^1(v)}.
$$
\end{conj}

We were not able to prove this conjecture but we have obtained the following result using an adequate Calderón-Zygmund decomposition that involves the $A_\infty$ constant of the weight.

\begin{teo} \label{teorema sawyer con Ap y Ainfty}
Let $v \in A_p$, $p \geq 1$, then there exists a dimensional constant $c$ such that
$$\left\| \frac{M_d (fv)}{v} \right\|_{L^{1, \infty}(v)}  \; \leq \; c\,[v]_{A_\infty}\max\{p,\,\log(e+[v]_{A_p})\} \|f\|_{L^1(v)}.$$
\end{teo}

%Taking into account that $[v]_{A_\infty} \leq c_n [v]_{A_p}$, $p \geq 1$, we derive the following corollary.

\begin{coro} \label{corolario para M y pesos Ap}
Let $v \in A_p$, $p \geq 1$, then there exists a dimensional constant $c$ such that
$$\left\| \frac{M_d (fv)}{v} \right\|_{L^{1, \infty}(v)}  \; \leq \; C_n\,[v]_{A_p}\max\{p,\,\log(e+[v]_{A_p})\} \|f\|_{L^1(v)}.$$
\end{coro}

We also try to improve the dependency on the weight constant using  some other refined  constants that were introduced in
\cite{hytonen-perez} and formalized in the work of Lerner and Moen \cite{lerner moen}.

\begin{teo} \label{teorema con constante mixta}
Let $v \in A_p$, $p \geq 1$, then there exists a dimensional constant $c$ such that
$$\left\| \frac{M_d (fv)}{v} \right\|_{L^{1, \infty}(v)}  \; \leq \; c\,p\;[v]_{(A_p)^{1/p}(A^{exp}_\infty)^{1/p'}}\log(e+[v]_{(A_p)^{1/p}(A^{exp}_\infty)^{1/p'}}) \|f\|_{L^1(v)}.$$
\end{teo}
 We remit to Section \ref{preliminares} for the definition of \,$[v]_{(A_p)^{1/p}(A^{exp}_\infty)^{1/p'}}$.

In this paper we will also study  similar problems for Calderón-Zygmund operators instead of the Hardy--Littlewood maximal function. In particular, we will improve the following theorem  from  \cite{hytonen-perez}.

\begin{teo} Suppose that  $T$ is a Calderón-Zygmund Operator, then there is a dimensional constant c such that for any $v\in A_1$
$$
\left\| \frac{T (fv)}{v} \right\|_{L^{1, \infty}(v)} \; \leq \; c\,[v]_{A_1}\log(e+[v]_{A_\infty})
\|f\|_{L^1(v)}.
% \int_{\R^n} |f(x)|\,dx.
$$
\end{teo}

This theorem improved the following result previously obtained in \cite{LOP2}.

$$
\left\| \frac{T (fv)}{v} \right\|_{L^{1, \infty}(v)} \; \leq \; c \,[v]_{A_1}\log(e+[v]_{A_1})
\|f\|_{L^1(v)}.
%\int_{\R^n} |f(x)|\,dx.
$$

 In section 5, we will give a version of Corollary \ref{corolario para M y pesos Ap} for Calderón-Zygmund operators. We will prove the following result.

\begin{teo}\label{teorema para un CZ}
Suppose that  $T$ is a Calderón-Zygmund Operator, then there is a dimensional constant c such that for any $v\in A_p$
$$\left\| \frac{T (fv)}{v} \right\|_{L^{1, \infty}(v)}  \; \leq \; c\,[v]_{A_p}\max\{p,\;\log(e+[v]_{A_p})\}
\|f\|_{L^1(v)}.$$
\end{teo}

\section{Preliminaries} \label{preliminares}

As usual a  weight will be a nonnegative locally integrable function. Given a weight $w$, $p\in (1,\infty)$ and a cube $Q$ we denote
$$
A_p(w;Q):=\left( \frac{1}{|Q|}\int_Q w\right)\left(\frac{1}{|Q|}\int_Q w^{1-p'}\right)^{p-1}= \frac{w(Q)\sigma(Q)^{p-1}}{|Q|^p},
$$
where $\sigma=w^{-\frac{1}{p-1}}$. When $p=1$ we define the limiting quantity as
$$
A_1(w;Q):= \left(\frac{1}{|Q|}\int_Q w\right)(\inf_Q w)^{-1}= \lim _{p\to 1}A_p(w,Q).
$$
For $p=\infty$ we will consider two constants. The first constant is defined as a limit of the  $A_p(w;Q)$ constants
$$
A_{\infty}^{exp}(w;Q):=\left( \frac{1}{|Q|}\int_Q w \right)\exp{\left( \frac{1}{|Q|}\int_Q \log w^{-1} \right)}=\lim_{p\to \infty}A_p(w,Q).
$$
To define the second constant we let
$$
A_{\infty}^{W}(w;Q):=\frac{1}{w(Q)}\int_Q M(\chi_Q w)
$$
and define:
$$
[w]_{A_p}= \sup_{Q }A_p(w;Q),
$$
$$
\|w\|_{A_\infty}=\sup_Q A_\infty^{exp}(w;Q)
$$
and
$$
[w]_{A_\infty}=\sup_Q A_{\infty}^{W}(w;Q).
$$
We write $w \in A_p$ if $[w]_{A_p}<\infty$ and $w \in A_\infty$ if $\|w\|_{A_\infty}<\infty$ or $[w]_{A_\infty}<\infty$. The constant $\|w\|_{A_\infty}$ was defined by Hru\v{s}\v{c}ev in \cite{Hruscev}. The constant $[w]_{A_\infty}$ was defined by Fujii in \cite{Fujii} and rediscovered by M. Wilson in \cite{Wilson87, Wilson08}, who also showed that both constants define the class $A_\infty$.
In \cite{hytonen-perez}, the authors proved the estimate
\begin{equation}
[w]_{A_\infty}\leq c_n\, \|w\|_{A_\infty}
\end{equation}
and provided examples showing that $\|w\|_{A_\infty}$ can be exponentially larger than $[w]_{A_\infty}$.

We now define the mixed type constants. Given $1\leq p < \infty$ and $\alpha$, $\beta\geq0$,  motivated by some results for the two weighted estimate for the maximal function in \cite{hytonen-perez}, Lerner and Moen in \cite{lerner moen}  defined  the following mixed constants:
\begin{equation*}
[w]_{(A_p)^\alpha(A_r)^\beta}= \sup_Q A_p(w;Q)^\alpha A_r(w;Q)^\beta,\;\;\; 1\leq r < \infty,
\end{equation*}
the exponential mixed constants:
\begin{equation}\label{mixedLM}
[w]_{(A_p)^\alpha(A_\infty^{exp})^{\beta}}= \sup_Q A_p(w;Q)^\alpha A_\infty^{exp}(w;Q)^\beta,
\end{equation}
and the Fujii-Wilson mixed constants:
$$
[w]_{(A_p)^\alpha(A_\infty^{W})^{\beta}}= \sup_Q A_p(w;Q)^\alpha A_\infty^{W}(w;Q)^\beta.
$$
%
%\end{defi}
If $\alpha>0$, the class of weights that satisfy
$$
[w]_{(A_p)^\alpha(A_\infty^{W})^{\beta}}<\infty,
$$
is simply the class  $A_p$, since
$$
\max([w]_{A_p}^{\alpha}, [w]_{A_\infty}^\beta)\leq [w]_{(A_p)^\alpha(A_\infty^{W})^{\beta}}\leq [w]_{A_p}^{\alpha+\beta}.
$$
Analogously, a weight $w$ satisfies $[w]_{(A_p)^\alpha(A_\infty^{exp})^{\beta}}<\infty$ if and only if $w$ is in $A_p$ such that  the inequality holds for the exponential mixed constant.
In \cite{lerner moen} the author showed that if $0< \alpha \leq \beta\leq 1$ and $w \in A_p$, then
\begin{equation}\label{desigualdadmixta}
[w]_{(A_p)^\alpha(A_\infty^{exp})^{1-\alpha}}\leq [w]_{(A_p)^\beta(A_\infty^{exp})^{1-\beta}}.
\end{equation}

We finish this section by defining the generalized $A_\infty$ class of weights $A_\infty (\mu)$ where $\mu$ is a doubling measure.
To do this we recall some well known definitions about generalized Hardy-Littlewood maximal operators. For a complete account, we refer the reader to \cite{duo, garcia curva-rubio de francia}.

Given a doubling measure $\mu$ we define the maximal operator $M_\mu$ by
$$
M_\mu f(x)=\sup_{Q \ni x}\frac{1}{\mu(Q)}\int_Q|f(y)|d\mu(y).
$$
For $1<p<\infty$, given a weight $w$ we say that $w \in A_p(\mu)$ if for all cubes $Q$,
$$
\left(\frac{1}{\mu(Q)}\int_Q w(x)\;d\mu(x)\right)\left(\frac{1}{\mu(Q)}\int_Q w(x)^{1-p'}\;d\mu(x)\right)^{p-1} \leq C.
$$
We say that  $w\in A_1(\mu)$ if
$$
M_\mu w(x)\leq Cw(x).
$$
 We denote the union of all the $A_p(\mu)$ classes by  $A_\infty (\mu)$, that is to say
 $$
 A_\infty (\mu)=\cup_{p\geq1}A_p(\mu).
 $$

Since $\mu$ is doubling, then $M_\mu$ is bounded on   $L^p(w d\mu)$, $1<p<\infty$, if and only if $w\in A_p(\mu)$. As usual when $\mu$ is the Lebesgue measure we omit the subscript  $\mu$  and write simply  $M$ or $A_p$. Also, if $\mu$ is absolutely continuous
given by the weight $u$ then we simply write $ A_p(u)$, $1\leq p\leq \infty$.
%Recall that if $w \in A_p$, then $w$ satisfies the reverse Hölder inequality
%$$
%\left(\frac{1}{|Q|}\int_Q w(x)^sdx\right)^{1/s}\leq \frac{C}{|Q|}\int_Q w(x)dx,
%$$
%%
%for some $s>1$ which depends only on $[w]_{A_p}$.
%
%\begin{defi} A weight  $w$  is a  $RH_s$ weight if
%%
%$$
%\left(\frac{1}{|Q|}\int_Q w(x)^sdx\right)^{1/s}\leq \frac{C}{|Q|}\int_Q w(x)dx,
%$$
%%
%and we denote the best constant in these inequalities by $[w]_{RH_s}$.
%
%\comment{\carlos{es necesario definir esta constante?   Donde se usa?}}
%
%Moreover, we denote by  $RH_\infty$ the class of weights such that
%%
%$$
%\frac{C}{|Q|}\int_Q w(x) dx \geq \esssup_Q w(x).
%$$
%%
%\end{defi}
%
%For every  $s>1$, $RH_\infty \subset RH_s$. For more information on these classes, see  \cite{C-UN}.

The next two  lemmas were proved in  \cite{C-UMaPe11}.
\begin{lema}\label{lema 2.1 de CUMP}
 If $u \in A_1$ and $v\in A_\infty(u)$, then  $uv\in A_\infty$. In particular, if   $v\in A_p(u)$,  $1\leq p<\infty$, then $uv\in A_p$.
\end{lema}

\begin{obs}\label{observacion 2.2 de CUMP}
  If $u\in A_1$, then $v \in A_\infty(u)$  if and only if $uv \in A_\infty$.
\end{obs}

%\begin{lema}\label{lema 2.3 de CUMP}
%  If $u\in A_1$ and $v\in A_p$, $1\leq p <\infty$, then there exists  $0<\epsilon_0<1$  depending only on  $[u]_{A_1}$ such that  $uv^{\epsilon}\in A_p$ for all $0<\epsilon<\epsilon_0$.
%\end{lema}
%The proof of following Lemma can see in \cite{C-UN}.
%
%\begin{lema}\label{lema 2.4 de CUMP}
%\
%\begin{enumerate}
%\item  $w\in A_\infty$ if and only if $w=w_1w_2$, where $w_1 \in A_1$ and $w_2 \in RH_\infty$.
%\item If $w\in A_1$, then $w^{-1}\in RH_\infty$.
%\item If $u,v \in RH_\infty$, then  $uv \in RH_\infty$.
%\end{enumerate}
%\end{lema}

%\comment{\carlos{No est\'an citados estos lemas de arriba}}

\begin{lema} \label{lema 2.5 de cump}
If $u \in A_1$ and $uv\in A_\infty$, then $v \in A_\infty$.
\end{lema}

\color{black}
\section{The $A_1$ case}

\begin{prueba}{Theorem \ref{teo u y v en A1 para la recta}}

As usual we denote $M^c$  the centered Hardy-Littlewood maximal operators and its
corresponding centered weighted $M_v^c$ maximal function.  Now, by standard arguments
$$\frac{M(fv)}{v} \approx \frac{M^c(fv)}{v} \leq \frac{M^cv}{v} M_v^c(f)\leq   \frac{Mv}{v} M_v^c(f)
\leq   [v]_{A_1} M_v^c(f)
$$
and then
$$ \left\|\frac{M(fv)}{v}\right\|_{L^{1,\infty}(v)}\leq
c_n [v]_{A_1} \,\left\|  M_v^c(f)  \right\|_{L^{1,\infty}(v)}
$$
$$
\leq
c_n [v]_{A_1} \,\left\|f  \right\|_{L^{1}(v)},
$$
by the Besicovitch covering lemma.

The proof will be completed  by showing that the linear exponent is the best possible. To see this, it is sufficient to considerer  $f(x)=\frac{1}{\delta}\chi_{(0,1)}(x)$ and
$v(x)= |x|^{\delta -1}$ where  $0<\delta<1$. Then standard computations shows that
$$[v]_{A_1}\sim \frac{1}{\delta}.$$
On the other hand, we can compute
 $$
M(fv)\geq \left\{
                     \begin{array}{ll}
                       \frac{1}{\delta}\frac{1}{x^{1-\delta}} & \text {if }  x \in (0,1)  \\ \\
                       \frac{1}{\delta^2}\frac{1}{x} & \text{if } x \in (1, \infty) \\ \\
                       \frac{1}{\delta^2}\frac{1}{1-x} & \text{if } x \in (-\infty,0) \\
                     \end{array}
                   \right.
$$
therefore $(0,\delta^{-2/\delta}) \subset \{x: M(fv)>v\}$. Continuing we have
$$v\{x: M(fv)>v\}\geqslant v(0, \delta^{-2/\delta})=\int_0^{\delta^{-2/\delta}} x^{\delta-1} dx = \frac{1}{\delta^3}=[v]_{A_1} \frac{1}{\delta^2},$$
 but $\int_{\R} f(x)v(x)dx=\int_0^1 \frac{1}{\delta}x^{\delta-1} dx =\frac{1}{\delta^2}$.

\end{prueba}

\begin{prueba}{Theorem \ref{teo que las constantes son como mimimo el producto}}
Let $f(x)=\frac{1}{\delta}\chi_{(0,1)}(x)$ and  define
$u(x)= \alpha \chi_{(0,1)}(x) + \chi_{(0,1)^c}(x)$, where  $0<\alpha <1$ and
$v(x)= |x|^{\delta -1}$, where $0<\delta<1$. Then standard computations shows that

$$[u]_{A_1}\sim \frac{1}{\alpha}\;\;\;\;\;\;\;\text{and}\;\;\;\;\;\;[v]_{A_1}\sim \frac{1}{\delta}.$$
Also, we have
 $$M(fv)\geq \left\{
                     \begin{array}{ll}
                       \frac{1}{\delta}\frac{1}{x^{1-\delta}} & \text {if }  x \in (0,1)  \\ \\
                       \frac{1}{\delta^2}\frac{1}{x} & \text{if } x \in (1, \infty) \\ \\
                       \frac{1}{\delta^2}\frac{1}{1-x} & \text{if } x \in (-\infty,0). \\
                     \end{array}
                   \right.
$$
Then, $(0,\delta^{-2/\delta}) \subset \{x: M(fv)>v\}$ and then
$$uv\{x: M(fv)>v\}\geqslant uv(1, \delta^{-2/\delta})=\int_1^{\delta^{-2/\delta}} x^{\delta-1} dx = \frac{1}{\delta}(\delta^{-2}-1)\thickapprox \frac{1}{\delta^3}.$$
On the other hand,
$$\int_\R f(x)u(x)v(x) dx= \frac{\alpha}{\delta}\int_0^1 x^{\delta-1} dx = \frac{\alpha}{\delta^2},$$
this proves $\varphi([u]_{A_1},[v]_{A_1})\gtrsim [u]_{A_1}[v]_{A_1}$.
\end{prueba}
 \begin{obs}
  When considering the case $\alpha = \delta$, we have  that $\varphi([u]_{A_1},[v]_{A_1})$ cannot be $\max([u]_{A_1},[v]_{A_1})$.
 \end{obs}

\section{The $A_p$ case} \label{Caso especial de la conjetura de Sawyer para $M_d$}

\begin{prueba}{Theorem \ref{teorema sawyer con Ap y Ainfty}}
Without loss of generality we may assume that  $f$ is nonnegative and bounded with compact support. Let $v\in A_p$ then $v\in A_r$, $r>p$ with $[v]_{A_r}\leq[v]_{A_p}$.  Fix $t>0$ and let $r>p$ be a parameter that will be chosen in a moment.  Since $v \in A_r$,  in particular, $vdx$  is a doubling weight.
Therefore, we can form the Calderón-Zygmund decomposition of $f$ at height $t$ with
respect to the measure $v(x)dx$. This yields a collection of disjoint dyadic maximal cubes $\{Q_j\}$, such that for all  $Q_j$:
$$
 t<\frac{1}{v(Q_j)}\int_{Q_j}f(x)v(x)dx \leq \frac{v(Q'_j)}{v(Q_j)v(Q'_j)}\int_{Q'_j}f(x)v(x)dx\leq 2^{nr}[v]_{A_r}t,
$$
where $Q'_j$ is the ancestor of $Q_j$ and where  the last inequality is obtained by using standard properties of the $A_p$ weights  (see Proposition 9.1.5 in \cite{Grafakos modern}) and by the maximality property of the $Q_j$.

Further,  if we let $\Omega:=\cup_j Q_j$, then  $f(x)\leq t$  for almost every  $x \in \R^n\setminus\Omega$. We decompose  $f$ as $g+b$, where

$$
 g(x)= \left\{
                     \begin{array}{ll}
                       \frac{1}{v(Q_j)}\int_{Q_j}f(x)v(x)dx & \text {if }  x \in Q_j  \\ \\

                        f(x)& \text{if } x \in  \R^n\setminus\Omega\\
                     \end{array}
                   \right.
$$

and let $b(x)=\sum_j b_j(x)$, with

$$
b_j(x)= \left(f(x)-{\frac{1}{v(Q_j)}}\int_{Q_j}f(x)v(x)dx\right) \chi_{Q_j}(x).
$$
If we used  this definitions, we have that $g(x)\leq 2^{nr}[v]_{A_r}t$ for almost every  $x \in \R^n$  and
$$
\int_{Q_j}b_j(x)v(x)dx=0.
$$

Following   \cite{cruz-uribe martell perez}  if $Q$ is a dyadic cube, then  $\forall x \in Q$,
$$
\frac{1}{|Q|} \int_Q f(x)v(x)dx= \frac{1}{|Q|}\int_Q g(x)v(x)dx +\frac{1}{|Q|}\int_Q b(x)v(x)dx \leq M_d(gv)(x)+ \widetilde{M}_d(bv)(x),
$$
where 
$$\widetilde{M}_d(h)(x)= \sup_{x\in Q}\left|\frac{1}{|Q|}\int_Q h(y)dy\right|.$$
Then, if the supremum is taken over all dyadic cubes containing $x$, we have
$$
M_d(fv)\leq M_d(gv)+ \widetilde{M}_d(bv).
$$
Now,
$$
v(\{x\in \R^n: \frac{M_d(fv)(x)}{v(x)}>t\})\leq v(\{x\in \R^n: \frac{M_d(gv)(x)}{v(x)}>t/2\})+
$$
$$
+v(\{x\in \Omega: \frac{\widetilde{M}_d(bv)(x)}{v(x)}>t/2\})+v(\{x\in \R^n\setminus \Omega : \frac{\widetilde{M}_d(bv)(x)}{v(x)}>t/2\})=I_1 +I_2+I_3.
$$
To estimate $I_1$ we will use the following improvement of Buckley's theorem (see \cite{buckley}) whose proof can be found in
\cite{hytonen-perez},
\begin{lema}
Let  $1<p<\infty$ and $v \in A_p$ then,
$$
\|M\|_{L^p(v)}\leq c_n p' \,     [v]_{A_p}^{1/p}\, [v^{1-p'}]_{A_{\infty} }^{1/p},
$$
where $c_n$ is a dimensional constant.
\end{lema}
We then have after applying Chebyshev inequality
$$
I_1\leq \frac{2^{r'}}{t^{r'}}\int_{\R^n} M(gv)^{r'} v^{1-r'}dx\leq \frac{c_n^{r'}}{t^{r'}}  r^{r'}[v]_{A_r}^{r'-1}[v]_{A_\infty} \int_{\R^n} g^{r'}v dx
.$$
 Since $g(x)\leq 2^{nr}[v]_{A_r}t$ and $[v]_{A_r} \leq [v]_{A_p}$, we have
$$
I_1\leq \frac{c_n^{r'}}{t}r^{r'}[v]_{A_\infty}[v]_{A_r}^{2r'-2}\int_{\R^n} g(x)v(x)dx\leq \frac{c_n^{r'}}{t}r^{r'}[v]_{A_\infty}[v]_{A_p}^{2r'-2}\int_{\R^n} g(x)v(x)dx
.$$
Finally, if we let
$$r= 1+\max\{ p,\,\log(e+[v]_{A_p})\},$$
 then
 $$r'=1+ \frac{1}{\max\{p,\log(e+[v]_{A_p})\}}.$$
and a computation shows that  $r^{r'}$  behaves like $\max\{p,\,\log(e+[v]_{A_p})\}$ and that $[v]_{A_p}^{2r'-2}$ is bounded.  Therefore,
$$
I_1 \leq \frac{C_n}{t}\,[v]_{A_\infty}\max\{p,\,\log(e+[v]_{A_p})\}\left(\int_{\R^n\setminus \Omega}f(x)v(x) dx+\sum_j \left(\frac{1}{v(Q_j)}\int_{Q_j}f(x)v(x)dx\right)v(Q_j)\right)
$$
$$
\leq \frac{C_n}{t}\,[v]_{A_\infty}\max\{p,\,\log(e+[v]_{A_p})\}\int_{\R^n}f(x)v(x)dx.
$$
The estimate for  $I_2$ follows immediately from the properties of the cubes  $Q_j$:
$$
I_2\leq v(\Omega)= \sum_j v(Q_j)\leq \sum_j \frac{1}{t}\int_{Q_j}f(x)v(x) dx\leq \frac{1}{t}\int_{\R^n}f(x)v(x) dx.
$$
Finally, we will prove that  $I_3=0$. To see this, fix  $x\in \R^n\setminus\Omega$, since  $b$ has support in  $\Omega$, to compute  $\widetilde{M}_d(bv)$  we only need to consider  cubes which intersect $\Omega$.
Fix such a cube $Q$, and for each $j$ either  $Q_j \subset Q$ or $Q\cap Q_j = \emptyset$. Then, since
$$\int_{Q_j}b_j(x)v(x)dx=0,$$
$$
\frac{1}{|Q|}\int_Q b(x)v(x)dx= \frac{1}{|Q|}\sum_j \int_{Q\cap Q_j} b_j(x)v(x)dx=\frac{1}{|Q|}\sum_{Q_j \subset Q} \int_{ Q_j} b_j(x)v(x)dx=0
.$$

\end{prueba}

We will use  the following lemma for mixed $A_p-A_\infty$ constants as defined in \eqref{mixedLM}.

\begin{lema}\label{acotar Ap con Ap Ainfty}
Let  $p>1$ and let $v \in A_p$, then
$$[v]_{(A_p)^{1/p}(A^{exp}_\infty)^{1/p'}}\leq [v]_{A_p} \leq [v]^p_{(A_p)^{1/p}(A^{exp}_\infty)^{1/p'}}.$$
\end{lema}

\begin{dem}

The second inequality  follows from  a simple consequence of a Jensen inequality:
$$
e^{\frac{1}{|Q|}\int_Q \log w(x) dx} \leq \frac{1}{|Q|}\int_Q w(x) dx,
$$
which implies $\left[(\frac{1}{|Q|}\int_Q w(x) dx)e^{\frac{1}{|Q|}\int_Q \log w^{-1}(x) dx}\right]^{p-1}\geq 1$ and then
$$
(\frac{1}{|Q|}\int_Q w(x) dx)(\frac{1}{|Q|}\int_Q w(x)^{1-p'} dx)^{p-1}\leq
$$
$$
\leq(\frac{1}{|Q|}\int_Q w(x) dx)^{p}(\frac{1}{|Q|}\int_Q w(x)^{1-p'} dx)^{p-1}(e^{\frac{1}{|Q|}\int_Q \log w^{-1}(x) dx})^{p-1},
$$
 whence we obtain
$$
[v]_{A_p} \leq [v]^p_{(A_p)^{1/p}(A^{exp}_\infty)^{1/p'}}.
$$

The first inequality also follows  from Jensen's inequality in the form
$$
e^{\frac{1}{|Q|}\int_Q \log w(x)^{-1} dx} \leq \left( \frac{1}{|Q|}\int_Q w(x)^{-\alpha} dx \right)^{1/\alpha} \qquad \alpha>0,
$$
considering the case $\alpha=p'-1$.

\end{dem}

We also need  the following lemma that will play an important role in the proof of Theorem \ref{teorema con constante mixta}. It is an improvement of Buckley's theorem (see \cite{buckley}) and the proof can be found in \cite{hytonen-perez}.

\begin{lema}\label{buckley mejorado para ApAinfty}
Let  $1<p<\infty$ and $v \in A_p$ then,
$$
\|M\|_{L^p(v)}\leq c_n p' \, [v^{1-p'}]_{(A_{p'})^{1/p'}(A_{\infty}^{exp})^{1/p}}
,$$
where $c_n$ is a dimensional constant.
\end{lema}

%The proof of Theorem \ref{teorema con constante mixta}
%similar a la del teorema \ref{teorema sawyer con Ap y Ainfty} solo nos detendremos en  cómo aplicamos los lemas antes mencionados y cómo perdemos la posibilidad de acotar con el máximo entre $p$ y la constante del peso y debemos quedarnos con el producto. Con lo cual, el teorema \ref{teorema sawyer con Ap y Ainfty} no se obtiene directamente de este.
%%
\color{black}

\begin{prueba}{Theorem \ref{teorema con constante mixta}}
The structure of the proof is the same as that of Theorem \ref{teorema sawyer con Ap y Ainfty}. The only difference is in the analysis of $I_1$. Indeed, combining Chebyshev inequality with Lemma   \ref{buckley mejorado para ApAinfty} we arrive to
$$
I_1\leq \frac{2^{r'}}{t^{r'}}\int_{\R^n} M(gv)^{r'} v^{1-r'}dx\leq \frac{2^{r'}}{t^{r'}} r^{r'}[v]_{(A_r)^{1/r}(A_{\infty}^{exp})^{1/r'}}^{r'}\int_{\R^n} g^{r'}v dx
$$
and since  $g(x)\leq 2^{nr}[v]_{A_r}t$ we have
$$
I_1\leq \frac{2^{r'(1+n)}}{t}r^{r'}[v]_{A_r}^{r'-1}[v]_{(A_r)^{1/r}(A_{\infty}^{exp})^{1/r'}}^{r'}\int_{\R^n} g(x)v(x)dx
.$$
As $r>p$, $[v]_{A_r} \leq [v]_{A_p}$ and by  \eqref{desigualdadmixta} $[v]_{(A_r)^{1/r}(A_{\infty}^{exp})^{1/r'}} \leq [v]_{(A_p)^{1/p}(A_{\infty}^{exp})^{1/p'}}$.  Finally, if we let
 $$r= 1+\max\{ p,\,\log(e+[v]_{A_p})\},$$
  then
$$r'=1+ \frac{1}{\max\{p,\log(e+[v]_{A_p})\}}.$$
It is easy to see that $r^{r'}$  behaves like
$$\max\{p,\,\log(e+[v]_{A_p})\},$$
that  $[v]_{(A_p)^{1/p}(A_{\infty}^{exp})^{1/p'}}^{r'}$ behaves like $[v]_{(A_p)^{1/p}(A_{\infty}^{exp})^{1/p'}}$ and that $[v]_{A_p}^{r'-1}$ is bounded by a universal constant. Moreover, since $2^{r'(1+n)}\leq 2^{2(1+n)}$ we have that
$$
I_1 \leq \frac{C_n}{t}\,[v]_{(A_p)^{1/p}(A_{\infty}^{exp})^{1/p'}}
\max\{p,\,\log(e+[v]_{A_p})\}\times
$$
$$
\times\left(\int_{\R^n\setminus \Omega}f(x)v(x) dx+\sum_j \left(\frac{1}{v(Q_j)}\int_{Q_j}f(x)v(x)dx\right)v(Q_j)\right)
.$$
Now by Lemma \ref{acotar Ap con Ap Ainfty} we have
$$
\max\{p,\,\log(e+[v]_{A_p})\}\leq \max\{p,\,p\log(e+[v]_{(A_p)^{1/p}(A_{\infty}^{exp})^{1/p'}})\}=p\log(e+[v]_{(A_p)^{1/p}(A_{\infty}^{exp})^{1/p'}})
$$
and then
$$
I_1 \leq\frac{C_n}{t}\,p\,[v]_{(A_p)^{1/p}(A_{\infty}^{exp})^{1/p'}}\log(e+[v]_{(A_p)^{1/p}(A_{\infty}^{exp})^{1/p'}})\int_{\R^n}f(x)v(x)dx.
$$
This concludes the proof of the theorem.

%La acotación de $I_2$ sale directamente de las propiedades de los $Q_j$:
%%
%$$
%I_2\leq v(\Omega)= \sum_j v(Q_j)\leq \sum_j \frac{1}{t}\int_{Q_j}f(x)v(x) dx\leq \frac{1}{t}\int_{\R^n}f(x)v(x) dx.
%$$
%%
%Finalmente veamos que $I_3=0$. Para ver esto, fijo $x\in \R^n\setminus\Omega$, como el soporte de $b$ esta en $\Omega$, para calcular $\widetilde M(bv)$ solo necesitamos considerar cubos cuya intersección con $\Omega$ sea no nula.
%Fijo tal $Q$, y para cada $j$, $Q_j \subset Q$ o $Q\cap Q_j = \emptyset$. Luego, como
% $$\int_{Q_j}b_j(x)v(x)dx=0,$$
%%
%$$
%\frac{1}{|Q|}\int_Q b(x)v(x)dx= \frac{1}{|Q|}\sum_j \int_{Q\cap Q_j} b_j(x)v(x)dx=\frac{1}{|Q|}\sum_{Q_j \subset Q} \int_{ Q_j} b_j(x)v(x)dx=0
%$$
%%

\end{prueba}

\section{Counterexamples} \label{contraejemplo}

In this section we show that inequality \eqref{casev=Mw}
%
%\begin{equation}\label{casev=Mw-bis}
%\left\|\frac{M(f)}{Mw}\right\|_{L^{1,\infty}(Mw)}\leq c\, \|f\|_{L^1(\R^n)},
%\end{equation}
is false. To do this we proceed by contradiction assuming that this inequality holds. We begin with the following  duality argument for any weight $w$,
$$
\left\| Tf\right\|_{ L^p(w)}= \sup_{h: \|h\|_{L^{p'}(w)}=1} |\int_{\mathbb{R}^{n}}Tf\,h\,wdx|
.$$
Fixing one of these $h$  we have
$$
\int_{\mathbb{R}^{n}}Tf\,h\,wdx= \int_{\mathbb{R}^{n}}f\, T^t(h\,w)dx=
\int_{\mathbb{R}^{n}}f\, \frac{T^t(h\,w)}{Mw}Mwdx
$$
and then
$$
|\int_{\mathbb{R}^{n}}Tf\,h\,wdx|
\leq \|f\|_{L^p(Mw)}\,  \left\| \frac{T^t(hw)}{Mw}  \right\|_{ L^{p'}(  Mw ) }=
\|f\|_{L^p(Mw)}\, \left\| T^tf\right\|_{ L^{p'}(Mw)^{1-p'})}
.$$

Now we will use the following lemma which is a particular version of the classical estimate of Coifman-Fefferman for any Calder\'on-Zygmund operator $T$: let $p\in (0,\infty)$ and let $w\in A_{\infty}$, then there is a constant
$c$ depending upon $p,T$ the $A_{\infty}$ constant of $w$ such that
\begin{equation}\label{CF}
\left\| Tf\right\|_{ L^p(w) } \leq\, c_{T,p,[w]_{A_{\infty}}}
\left\| Mf\right\|_{ L^p(w) }.
\end{equation}
Then as a consequence we have the special situation: Let \, $w$ \, be {\it any weight} and let \quad
$p\in (1,\infty)$. \,Then, there is a constant depending only on $p$ and $T$ such that:
\begin{equation}\label{specialCF}
\left\| Tf\right\|_{ L^p(Mw)^{1-p})} \le\, c\,
\left\|Mf\right\|_{  L^p(Mw)^{1-p})}.
\end{equation}
This follows from \eqref{CF} and the fact $(Mw)^{1-p}\in A_{\infty}$. Indeed, since $(Mw)^{1-p}=(Mw)^{ \delta(1-2p)  }\in A_{2p}$  and since  $\delta=\frac{p-1}{2p-1}<\frac12$ we have $[(Mw)^{1-p}]_{A_{\infty}} \leq [(Mw)^{\delta}]^{2p-1}_{A_1} \leq c_n^{p}$.

It should be mentioned that  \eqref{specialCF} was improved in   \cite{LOP3} and later  in \cite{LOP1}. In these papers the relevance was the sharpness of the constant $c$ in term of $p$ which behaves linearly in $p$, but is not important in our context. See also \cite{Re 2013} for a similar estimate within the fractional integrals context.

Then, since $T^t$ is also a Calder\'on--Zygmund operator we apply \eqref{specialCF}
$$
|\int_{\mathbb{R}^{n}}Tf\,h\,wdx|
\leq
c_{p,T}\,
\|f\|_{L^p(Mw)}\, \left\| M(hw)\right\|_{ L^{p'}(Mw)^{1-p'})}
=
c_{p,T}\,
\|f\|_{L^p(Mw)}\,  \left\| \frac{M(hw)}{Mw}  \right\|_{ L^{p'}(  Mw ) }
.$$
We now apply  \eqref{casev=Mw} which is equivalent to
\begin{equation*}
\left\|\frac{M(fw)}{Mw}\right\|_{L^{1,\infty}(Mw)}\leq c\, \|f\|_{L^1(w)},
\end{equation*}
and since the operator $f\to \frac{M(fw)}{Mw}$ is trivially bounded on $L^{\infty}$ with constant $1$  we apply the Marcinkiewicz's interpolation theorem to deduce
$$
|\int_{\mathbb{R}^{n}}Tf\,h\,wdx|
\leq c_p \|f\|_{L^p(Mw)}\,\|h\|_{L^{p'}(w)}   = c_p \|f\|_{L^p(Mw)}.
$$
Then, for any Calder\'on--Zygmund operator $T$ and arbitrary weight $w$ we have produced the estimate
$$
\left\| Tf\right\|_{ L^p(w)} \leq c_p \|f\|_{L^p(Mw)}.
$$
However, this inequality is well known to be false for any $p\in(1,\infty)$  as was shown by M. Wilson in
\cite{Wilson89}  for the simplest case, namely the Hilbert transform.

\color{black}

\section{Calderón-Zygmund integral operator}

In this section, we will show the following inequality

$$
\left\| \frac{T (fv)}{v} \right\|_{L^{1, \infty}(v)}  \; \leq \; C_n\,[v]_{A_p}\max\{p,\;\log(e+[v]_{A_p})\} \int_{\R^n} |f(x)|v(x)\,dx.
$$
For the proof of this inequality are need  the following two results. The first result was proved  in \cite{Hytonen} and the second result can be found in \cite[p.413]{garcia curva-rubio de francia}.

\begin{teo} \label{prueba de la conjetura A2}
 Let $1<p<\infty$,   $w$ an $A_p$-weight  and $T$ is a  Calderón-Zygmund operator, then
$$
\|Tf\|_{L^p(w)} \leq c_n \;p \;p'\;[w]_{A_p}^{\max(1,\frac{1}{p-1})}
.$$
\end{teo}
\begin{lema}\label{lema para acotar la parte mala de GCyRF}
Let  $w$ be a weight. There is a dimensional constant $c_d$ such that for all cube  $Q$ and for all function  $f$ supported in a cube $Q$ with $\int_Q f(x)dx=0$, the following inequality is holds:
$$
\int_{\R^n\setminus2Q}|Tf(y)|w(y)dy \leq c_n \int_{Q}|f(y)|Mw(y)dy
.$$
\end{lema}
The structure of the proof of the  theorem \ref{teorema para un CZ} is similar to that of Theorem \ref{teorema sawyer con Ap y Ainfty}.

\begin{prueba}{Theorem \ref{teorema para un CZ}}
Without loss of generality we will assume that $f$ is bounded and has compact support.
 Since $v\in A_p$, then $\forall r>p$, we have $v\in A_r$, with $[v]_{A_r}\leq[v]_{A_p}$.\\
 Fix $t>0$. For now let  $r>p$ be arbitrary, we will assign a specific value to $r$.  Since $v \in A_r$,  in particular, $vdx$  is a doubling weight.
Therefore, we can form the Calderón-Zygmund decomposition of $f$ at height $t$ with
respect to the measure $vdx$. This yields a collection of disjoint dyadic maximal cubes $\{Q_j\}$, such that for all  $Q_j$:
$$
t<\frac{1}{v(Q_j)}\int_{Q_j}f(x)v(x)dx \leq \frac{v(Q'_j)}{v(Q_j)v(Q'_j)}\int_{Q'_j}f(x)v(x)dx\leq 2^{nr}[v]_{A_r}t,
$$
where as before $Q'_j$ is the ancestor of $Q_j$ and where  the last inequality is obtained by using standard properties of the $A_p$ weights  (see Proposition 9.1.5 in \cite{Grafakos modern}) and by maximal property of the $Q_j$. Further,  if we let $\Omega:=\cup_j Q_j$, then  $f(x)\leq t$  for almost every  $x \in \R^n\setminus\Omega$. We decompose  $f$ as $g+b$, where

$$
 g(x)= \left\{
                     \begin{array}{ll}
                       \frac{1}{v(Q_j)}\int_{Q_j}f(x)v(x)dx & \text {if }  x \in Q_j  \\ \\

                        f(x)& \text{if } x \in  \R^n\setminus\Omega\\
                     \end{array}
                   \right.
$$

and let $b(x)=\sum_j b_j(x)$, with

$$
b_j(x)= \left(f(x)-{\frac{1}{v(Q_j)}}\int_{Q_j}f(x)v(x)dx\right) \chi_{Q_j}(x).
$$
If we used  this definitions, we have that $g(x)\leq 2^{nr}[v]_{A_r}t$ for almost every  $x \in \R^n$  and
$$
\int_{Q_j}b_j(x)v(x)dx=0
.$$

%Fijamos\, $t>0$ y sea $r>p$, al cual después le daremos un valor especifico, como $v \in A_r$, en particular $v(x)dx$ es una medida duplicante y puedo hacer una descomposición de Calderón-Zygmund de $f$ a nivel $t$ respecto a la medida $v(x)dx$.
%
%Encontramos así una colección de cubos diádicos maximales $\{Q_j\}$, tales que para todo $Q_j$:
%%
%$$
%t<\frac{1}{v(Q_j)}\int_{Q_j}f(x)v(x)dx \leq \frac{v(2Q_j)}{v(Q_j)v(2Q_j)}\int_{2Q_j}f(x)v(x)dx\leq 2^{nr}[v]_{A_r}t
%$$
%%
%donde la última desigualdad se obtiene por la proposición 9.1.5 de \cite{Grafakos modern} y por la maximalidad de los $Q_j$.
%
%Si llamamos $\Omega:=\cup_j Q_j$ y $\tilde{\Omega}:=\cup_j 2Q_j$, entonces $f(x)\leq t$ en casi todo punto $x \in \R^n\setminus\Omega$. Si descomponemos $f$ como $g+b$, donde
%
%%
%$$
% g(x)= \left\{
%                     \begin{array}{ll}
%                       \frac{1}{v(Q_j)}\int_{Q_j}f(x)v(x)dx & \text {Si }  x \in Q_j  \\ \\
%
%                        f(x)& \text{Si } x \in  \R^n\setminus\Omega\\
%                     \end{array}
%                   \right.
%$$
%%
%
%y sea $b(x)=\sum_j b_j(x)$, con
%
%$$
%b_j(x)= \left(f(x)-{\frac{1}{v(Q_j)}}\int_{Q_j}f(x)v(x)dx\right) \chi_{Q_j}(x).
%$$
%Usando estas definiciones, tenemos que $g(x)\leq 2^{nr}[v]_{A_r}t$ en casi todo punto $x \in \R^n$ y que
%$$\int_{Q_j}b_j(x)v(x)dx=0.$$

Then, since $T$ is a sublineal operator we have that
$$
v(\{x\in \R^n: \frac{|T(fv)(x)|}{v(x)}>t\})\leq v(\{x\in \R^n: \frac{|T(gv)(x)|}{v(x)}>t/2\})+
$$
$$
+v(\{x\in \tilde{\Omega}: \frac{|T(bv)(x)|}{v(x)}>t/2\})+v(\{x\in \R^n\setminus \tilde{\Omega} : \frac{|T(bv)(x)|}{v(x)}>t/2\})=I_1 +I_2+I_3
.$$
If first we used Chebyshev inequality and later we apply Theorem \ref{prueba de la conjetura A2} bearing in mind that as $v \in A_r$ we have $v^{1-r'} \in A_{r'}$ with
  $$[v^{1-r'}]_{A_{r'}}=[v]_{A_r}^{r'-1}.$$
 Since the exponent of the constant $[v]_{A_{r'}}$ in Lemma \ref{prueba de la conjetura A2} is different if  $p>2$ or $p\leq 2$, we have divided the proof into two cases.

 \textbf{Case $p>2$:}
 In this case, as $r>2$, we have $r'< 2$ and $\max(1,\frac{1}{r'-1})=\frac{1}{r'-1}$.
$$
I_1\leq \frac{2^{r'}}{t^{r'}}\int_{\R^n} |T(gv)(x)|^{r'} v(x)^{1-r'}dx\leq \frac{c_n^{r'}}{t^{r'}} r^{r'}[v]_{A_r}^{r'}\int_{\R^n} g(x)^{r'}v(x) dx,
$$
Then, since $g(x)\leq 2^{nr}[v]_{A_r}t$  and $[v]_{A_r}\leq[v]_{A_p}$ we obtened that
$$
I_1\leq \frac{c_n^{r'}}{t}r^{r'}[v]_{A_r}^{2r'-1}\int_{\R^n} g(x)v(x)dx\leq \frac{2^{r'(1+n)}}{t}r^{r'}[v]_{A_p}^{2r'-1}\int_{\R^n} g(x)v(x)dx.
$$
As $r>p>2$, we choose
$$r= 1+\max\{ p,\,\log(e+[v]_{A_p})\},$$
 then
  $$2>r'=1+ \frac{1}{\max\{p,\log(e+[v]_{A_p})\}}.$$
   For this reason, $r^{r'}$ behaves like $\max\{p,\,\log(e+[v]_{A_p})\}$ and $[v]_{A_p}^{2r'-1}$ like $[v]_{A_p}$.
$$
I_1 \leq \frac{C_n}{t}\,[v]_{A_p}\max\{p,\,\log(e+[v]_{A_p})\}\left(\int_{\R^n\setminus \Omega}f(x)v(x) dx+\sum_j \left(\frac{1}{v(Q_j)}\int_{Q_j}f(x)v(x)dx\right)v(Q_j)\right)
$$
$$
\leq \frac{C_n}{t}\,[v]_{A_p}\max\{p,\,\log(e+[v]_{A_p})\}\int_{\R^n}f(x)v(x)dx.
$$

 \textbf{Case $p\leq2$:}
 We choose $r=1+2\log(e+[v]_{A_p})>2\geq p$,  thus $$r'=1+\frac{1}{2\log(e+[v]_{A_p})}< 2$$
   and $\max(1,\frac{1}{r'-1})=\frac{1}{r'-1}$. We can now proceed analogously to the previous case,
 $$
 I_1\leq \frac{c_d^{r'}}{t}r^{r'} [v]_{A_p}^{2r'-1}\int_{\R^n} g(x)v(x)dx,
 $$
 therefore
 $$
 I_1\leq \frac{c_d}{t}\,\log(e+[v]_{A_p})\int_{\R^n} f(x)v(x)dx.
 $$

The estimate for  $I_2$ follows immediately from the properties of the cubes  $Q_j$ and from the following inequality
 $$v(2Q)\leq 2^{np}[v]_{A_p}v(Q).$$
$$
I_2\leq v(\tilde{\Omega})\leq \sum_j v(2Q_j)\leq   2^{np}[v]_{A_p}\sum_j \frac{1}{t}\int_{Q_j}f(x)v(x) dx\leq 2^{np}[v]_{A_p}\frac{1}{t}\int_{\R^n}f(x)v(x) dx.
$$
Finally, to be able to  estimate  $I_3$ we used  Lemma \ref{lema para acotar la parte mala de GCyRF} with $w\equiv1$.
$$
I_3\leq \frac{2}{t}\int_{\R^n\setminus \tilde{\Omega}}|T(bv)(x)|dx\leq \frac{2}{t}\sum_j\int_{\R^n\setminus 2Q_j}|T(b_jv)(x)|dx\leq \frac{2}{t}\sum_j\int_{Q_j}|b_j(x)|v(x)dx,
$$
 if we used $b_j$'s definition, we have that
$$
I_3\leq \frac{c_d}{t}\|fv\|_{L^1(\R^n)}.
$$
\end{prueba}

\section{An adaptation of Sawyer's proof with control of the constant}\label{dependencia de las constantes u y v}
In this appendix, we will prove Theorem \ref{terorema de costantes cuadrado y a la cuarta} using a method similar as the one considered  in \cite{cruz-uribe martell perez} for proof of  theorem 1.4.

The statement of the theorem assumes that the weights  belong to the $A_1$ class of weights. This weights satisfy
a reverse H\"older inequality, namely if $w\in A_1$, then there are two constants
$r,c>1$  such that
$$ \left(\frac{1}{|Q|}\int_Q w^r\right)^{1/r} \le \frac{c}{|Q|}\int_Q
w.$$
However, in the classical proofs there is a bad dependence on the constant
$c=c(r,[w]_{A_1})$ and we need a more precise
estimate to get our results.

\begin{lema}\label{RH} Let $w\in A_1$, and let
$r_{w}=1+\frac{1}{2^{n+1}[w]_{A_1}}$. Then for any cube $Q$
$$ \left(\frac{1}{|Q|}\int_Q w^{r_w}\right)^{1/r_w} \le \frac{2}{|Q|}\int_Q
w.$$
As a consequence we have that for any cube $Q$ and for any measurable set $E\subset Q$
\begin{equation*}
 \frac{w(E)}{w(Q)} \leq 2 \left(\frac{|E|}{|I|}\right)^{\epsilon_w },
 \end{equation*}
where $\epsilon_w =\frac1{1+2^{n+1}[w]_{A_1} }. $
\end{lema}

The proof of this reverse H\"older inequality can be found in \cite{LOP1} and the consequence is an application of H\"older's inequality.

\begin{prueba}{Theorem \ref{terorema de costantes cuadrado y a la cuarta}}

%\begin{dem}
Fix $t>0$ and define $g = fv/t$. Then, it is sufficient to show that
\begin{equation}
 uv(\{ x \in {\R}^n  :\; M_d (g)(x) > v(x) \}) \; \leq \; C \int_{\R^n} |g(x)|u(x)\,dx, \label{b1}
\end{equation}
for any function  $g$ bounded with compact support.

Fix $a> 2^n$. For each  $k \in \Z$, let $\{ I_j^k\}$ be the collection of maximal, disjoint dyadic cubes whose union is the set
$$\Omega_k = \{x\in\R^n:M_d v(x)>a^k\}\cap \{x\in\R^n:M_dg(x)>a^k\}.$$
This decomposition exists since  $g$ is bounded and has compact support, so the second set is contained in the union of maximal dyadic cubes. Define

$$\Gamma= \{ (k,j): \,|I_j^k \cap\{x: v(x)\leq a^{k+1}\}|\,>\,0\}.$$

 As $v \in A_1$, we have $Mv(x)\leq [v]_{A_1}v(x)$ almost everywhere. Hence, for $(k,j) \in \Gamma$
\begin{equation}
\frac{a^k}{[v]_{A_1}} \leq \frac{1}{[v]_{A_1}} \essinf_{x \in I_j^k} M_d v(x)\leq \essinf_{x \in I_j^k}v(x)\leq \frac{1}{|I^k_j|}\int_{I^k_j}v(x)dx\leq [v]_{A_1} a^{k+1}. \label{b2}
\end{equation}
(Intuitively, if $(k,j)\in\Gamma$, then  $I_j^k$  behaves like a cube from the   Calderón-Zygmund decomposition of $v$ at height $a^k$). Then up to a set of measure zero we have the following inclusions: for each $k$,
$$
\{x\in\R^n: a^k<v(x)\leq a^{k+1}\}\,\cap\,\{x\in\R^n: M_d g(x)>v(x)\}\,\subset\,\bigcup_{j:(k,j)\in\Gamma} I_j^k
.$$
Combining this with (\ref{b2}) we get that
$$uv(\{x\in\R^n :M_d g(x)>v(x)\})\leq a [v]_{A_1} \sum_{(k,j)\in\Gamma} {|I_j^k|}^{-1}v(I_j^k)u(I_j^k).$$
Fix $N<0$ and define  $\Gamma_N =\{(k,j) \in \Gamma : k\geq N \}$.  We will show that
$$
\sum_{(k,j)\in \gamma_N}{|I_j^k|}^{-1}v(I_j^k)u(I_j^k)\,\leq \; C \int_{\R^n} |g(x)|u(x)dx
.$$
where the constant $C$ does not depended of $N$. Inequality (\ref{b1}) then follows if we take the limit as $N\rightarrow -\infty$.
To prove this, we are going to replace the set of cubes $\{I_j^k\}$ by a subset with better properties. First, since   $v \in A_1$
we can apply Lemma \ref{RH} and there exists  $\epsilon = {(1+2^{n+1}[v]_{A_1})}^{-1} >0$ such that given any cube $I$
 and $E\subset I$,
 \begin{equation}
 \frac{v(E)}{v(I)} \leq 2 \left(\frac{|E|}{|I|}\right)^\epsilon. \label{b3}
 \end{equation}
 Fix $\delta$ such that $0<\delta<\epsilon$.  Define $\Delta_N =\{I_j^k : (k,j) \in \Gamma_N\}$. The cubes in  $\Delta_N$ are all dyadic, so they are either paiwise disjoint or ine is contained in the other. For $k>t$, since $\Omega_k \subset \Omega_t$ and since the cubes  $I_j^k$ are maximal in  $\Omega_k$, if $I_s^t \cap I_j^k \neq \emptyset$, then  $I_j^k \subset I_s^t$. In particular, each cube  $I_j^k \in \Delta_N$ is contained in  $\cup_j I_j^N \subset \{x: M_d g(x)>a^N\}$. As we noted above, the last set is bounded, so  $\Delta_N$ contains a maximal disjoint subcollection of cubes.

We form a sequence of sets  $\{G_n\}$ by induction. Let  $G_0$ be the set of all pairs $(k,j)\in \Gamma_N$ such that $I_j^k$ is maximal in $\Delta_N$.  For $n\geq0$, given the set   $G_n$, define the set  $G_{n+1}$ to be the set of pairs   $(k,j) \in \Gamma_N$ such that there exists  $(t,s) \in G_n$ with $I_j^k \subsetneq I_s^t$ and \\
 \begin{equation}
 \frac{1}{|I_j^k|} \int_{I_j^k} u(x)dx\,>\, a^{(k-t)\delta} \frac{1}{|I_s^t|}\int_{I_s^t}u(x)dx, \label{b4}
 \end{equation}
 \begin{equation}
 \frac{1}{|I_i^l|} \int_{I_i^l} u(x)dx \leq a^{(l-t)\delta} \frac{1}{|I_s^t|}\int_{I_s^t}u(x)dx. \label{b5}
 \end{equation}
Whenever $(l,i) \in \Gamma_N$ and $I_j^k \subsetneq I_i^l \subset I_s^t$.

 Let $P=\cup_{n\geq 0}G_n$. Given $(s,t) \in P$, we refer to the cube  $I_s^t$ as a principal cube. Since every cube in $\Delta_N$ is contained in a maximal cube, every cube in $\Delta_N$  is contained in one or more principal cubes.

 To continue, we divide  the proof into several steps the same form.  We will only  look at the behavior of the $A_1$-constants and we give the main ideas of the steps.

 \textbf{Step 1}\\
 We claim that
 \begin{equation}
 \sum_{(k,j)\in \Gamma_N} {|I_j^k|}^{-1}v(I_j^k)u(I_j^k)\,\leq \, C_{\epsilon} \sum_{(k,j)\in P}{|I_j^k|}^{-1}v(I_j^k)u(I_j^k). \label{b6}
 \end{equation}
 To prove this. Fix $(t,s) \in P$ and let $Q=Q(t,s)$ be the set of indices $(k,j) \in \Gamma_N$ such that $I_j^k \subset I_s^t$ and  $I_s^t$ is the smallest principal cube containing  $I_j^k$. In particular, each  $I_j^k$ is not a principal cube unless it equals $I_s^t$.

 So by  (\ref{b5}) and since  $I_j^k \subset \{x:M_d v(x)> a^k\}$,
 $$\sum_{(k,j)\in Q} {|I_j^k|}^{-1}v(I_j^k)u(I_j^k) \,\leq \; {|I_s^t|}^{-1} u(I_s^t) \sum_{k\geq t} a^{(k-t)\delta} v( I_s^t \cap\{x: M_d v(x)>a^k\}).$$
 By (\ref{b3}), (\ref{b2}), and since  $v\in A_1$,
 $$v( I_s^t \cap\{x: M_d v(x)>a^k\})\,\leq \; 2 [v]_{A_1}^{2\epsilon}a^{\epsilon}a^{(t-k)\epsilon}v(I_s^t).$$
 Combining these inequalities, we see that
 $$\sum_{(k,j)\in Q} {|I_j^k|}^{-1}v(I_j^k)u(I_j^k) \,\leq \; C_{\epsilon} {|I_s^t|}^{-1} u(I_s^t)v(I_s^t),$$
  where $C_{\epsilon} = \frac{a^{2\epsilon -\delta}}{(a^{\epsilon-\delta}-1)}2[v]_{A_1}^{2\epsilon}$. \\
 If we now sum over all  $(s,t)\in P$, we get (\ref{b6}) since $\cup_{(t,s)\in P}Q(t,s)=\Gamma_N$.\\

 \textbf{Step 2}
  For each $k$, let $\{J_i^k\}$ be the collection of maximal disjoint cubes whose union is  $\{x:M_d g(x)>a^k\}$. Then
 $$a^k < \frac{1}{|J_i^k|} \int_{J_i^k} g(x)dx.$$
  For each $j$, $I_j^k \subset \{x:M_d g(x)>a^k\}$, so there exists a unique  $i=i(j,k)$ such that $I_j^k \subset J_i^k$.
 Hereafter, the index $i$ will always be this function of $(k,j)$. Hence, by (\ref{b6}) and by (\ref{b2}),
 $$\sum_{(k,j)\in\Gamma_N} {|I_j^k|}^{-1}v(I_j^k)u(I_j^k)\,\leq\,C_{\epsilon} a [v]_{A_1}\int_{\R^n} h(x)g(x)dx,$$
 where $h(x)=\sum_{(k,j)\in P}{|J_i^k|}^{-1}\chi_{J_i^k}(x)u(I_j^k)$.

 To complete the proof we will show that for each  $x$, $h(x)\leq C u(x)$.
 Fix $x \in \R^n$; without loss of generality we may assume that  $u(x)$ is finite. For each  $k$, there exists al most one cube  $J_b^k$ such that  $x\in J_b^k$. If it exist, denote this cube by  $J^k$.

 Define $P_k=\{(k,j)\in P: I_j^k \subset J^k\}$, and $G=\{k: P_k\neq\emptyset\}$. We form a sequence  $\{k_m\}$ by induction. If $k\in G$, then $k\geq N$, so let $k_0$ be the least integer in  $G$. Given $K_m$, $m\geq0$, choose $k_{m+1}>k_m$ in  $G$ such that
 \begin{equation}
 \frac{1}{|J^{k_{m+1}}|} \int_{J^{k_{m+1}}}u(y)dy > \frac{2}{|J^{k_m}|} \int_{J^{k_m}}u(y)dy, \;\;\;\;\;\;\label{b7}
 \end{equation}
 \begin{equation}
 \frac{1}{|J^l|}\int_{J^l} u(y)dy \leq \frac{2}{|J^{k_m}|}\int_{J^{k_m}} u(y)dy,\;\;\;\;\;k_m\leq l<k_{m+1}, l\in G. \label{b8}
 \end{equation}
 Since  $u(x)$ is finite, the sequence  $\{K_m\}$ only contains a finite number of terms. Then by (\ref{b8}), we have

 $$h(x)\leq \sum_m \frac{2}{|J^{k_m}|}\int_{J^{k_m}} u(y)dy\sum_{l\in G, k_m\leq l<k_{m+1}} \sum_{(l,j)\in P_l} \frac{u(I_j^l)}{u(J^l)}.$$
 we claim that
 \begin{equation}
 \sum_{l\in G, k_m\leq l<k_{m+1}} \sum_{(l,j)\in P_l} \frac{u(I_j^l)}{u(J^l)}\leq C_9, \label{b9}
 \end{equation}
given this, we would be done: since the sequence  $\{k_m\}$ is finite, let $m$ be the largest index. Then by  (\ref{b7}) and (\ref{b9}),
$$h(x)\leq 2C_9 (2-{(\frac{1}{2}})^m){[u]}_{A_1} u(x).$$
Therefore, to complete the proof we must show  (\ref{b9}). We do this in two steps.\\

\textbf{Step 3}
They proved  in \cite{cruz-uribe martell perez} that if  $(l,j)\in P_l$, $k_m\leq l<k_{m+1}$, then
\begin{equation}
\frac{1}{|I_j^l|} \int_{I_j^l} u(y)dy \,>\,\frac{a^{(l-k_m)\delta}}{2{[u]_{A_1}}} \frac{1}{|J^l|}\int_{J^l} u(y) dy. \label{b10}
\end{equation}

\textbf{Step 4}
We will now  (\ref{b9}). By (\ref{b10}) and again since  $u\in A_1$, if $y\in I_j^l$, Then
$$\lambda=\frac{a^{(l-k_m)\delta}}{2{[u]}_{A_1}} \frac{u(J^l)}{|J^l|} \frac{1}{[u]_{A_1}} \,<\, u(y);$$
hence,
$$\cup_{j:(l,j)\in P_l} I_j^l \subset \{x\in J^l: u(x)>\lambda\}.$$
For $l$ fixed the cubes $I_j^l$ are disjoint. Therefore, since  $u\in A_1$ there exist  $\nu= (1+2^{n+1}[u]_{A_1})^{-1}$ such that
$$\sum_{j:(l,j)\in P_l}u(I_j^l)\leq 2^{1+\nu}u(J^l)[u]_{A_1}^{2\nu} a^{(k_m -l)\delta\nu}.$$
Therefore, we have that
$$\sum_{l\in G, k_m\leq l<k_{m+1}} \sum_{(l,j)\in P_l} \frac{u(I_j^l)}{u(J^l)}\leq C_9,$$
where $C_9= 2^{1+\nu}[u]_{A_1}^{2\nu} \frac{a^{\delta\nu}}{a^{\delta\nu}-1}$. Then, the constant $C$ of the theorem  \ref{terorema de costantes cuadrado y a la cuarta} behaves like
$$\frac{a^{\epsilon-\delta}}{a^{\epsilon-\delta}-1}2^{3+\nu}a^{\epsilon+2}[v]_{A_1}^{2\epsilon+2}[u]_{A_1}^{2\nu+1}\frac{a^{\delta\nu}}{a^{\delta\nu}-1} (2-(\frac{1}{2})^m) \approx [v]_{A_1}^4 [u]_{A_1}^2.$$

\end{prueba}

\section{Acknowledgement}

The second author is supported by the Spanish government grant MTM-2014-53850-P, 
the first and  third authors are  supported by Universidad Nacional del Sur and CONICET.

Finally, the authors would like to thank the referee whose suggestions and comments 
have been very helpful to improve the presentation of the paper.

\small

\markright{}

\emph{Sheldy Ombrosi, Departamento de Matem\'aticas, Universidad Nacional del Sur, 8000 Bahia Blanca, Argentina.}
e-mail address: sombrosi@uns.edu.ar

\emph{Carlos P\'erez, Department of Mathematics, University of the Basque Country and
Ikerbasque,  Bilbao, Spain}
e-mail address: carlos.perezmo@ehu.es

\emph{Jorgelina Recchi, Departamento De Matem\'aticas, Universidad Nacional Del Sur, 8000 Bahia Blanca, Argentina.}
e-mail address: drecchi@uns.edu.ar

\end{document}